\documentclass[12pt,a4paper]{article}
\usepackage{amsfonts}
\usepackage{amssymb}
\usepackage{amsmath}
\usepackage[english,francais]{babel}
\usepackage{icomma}
\usepackage[T1]{fontenc}
\usepackage{sectsty}
\usepackage{lipsum}
\allsectionsfont{\centering}
\usepackage{theorem}
\everymath{\displaystyle} 
\newtheorem{Theo}{Theorem}[section]

\newtheorem{Lem}{Lemma}[section]

\newtheorem{Rem}{Remark}[section]

\title{Distribution laws of smooth divisors} \,
\author{S. Nyandwi \;and \; A. Smati\;}

\date{ }
\begin{document}
\maketitle

\selectlanguage{english}

\begin{abstract}
A classical result due to Deshouillers, Dress and Tenenbaum asserts that on average the distribution of the divisors of the integers follows  the arcsine law.
In this paper, we investigate the distribution of  smooth divisors of the integers,  that is, those  divisors  which are free of large prime factors. 
We show that on average these divisors are distributed according to a probability law that we will describe.
\end{abstract}
\section{Introduction}
Let $n\geqslant1$ be an interger. We denote by $P(n)$  the largest prime divisor of $n\geqslant 2$ and we set $P(1)=1$. Let  $y\in ]1,+\infty[$ be a real number. Consider the set of   $y$-smooth divisors of $n$, that is, those  divisors of  $n$ which are  free of prime factors  exeeding $y$.
 $$\mathcal{D}_{n,y}:=\big\{d| n : P(d)\leqslant y\big\},$$
 and denote by $\tau(n,y)$ its cardinality. For each integer $n\geqslant1$ and for each real number $y>1$, we define the random variable $$X_{n,y}:\mathcal{D}_{n,y}\longrightarrow [0,1],$$
 which takes the values $\log d / \log n$ with uniform probability  $1/ \tau(n,y)$ and  for  $v\in[0,1]$, its distribution function
 $$F_{n,y}(v):=\mathbb{P}(X_{n,y}\leqslant v)=\frac{1}{\tau(n,y)}\,\sum_{d|n, d\leqslant n^{v},P(n)\leqslant y} 1.$$
It is easy to see that the sequence $(F_{n,y})_{n\geqslant1}$ does not converge pointwise  in $[0,1]$. We consider its mean in the interval $[0,1]$
 \begin{eqnarray}
 \frac{1}{x}\sum_{n\leqslant x}F_{n,y}(v)=\frac{1}{x}\sum_{n\leqslant x}\mathbb{P}(X_{n,y}\leqslant v)\cdot
 \end{eqnarray}
 The aim of this  paper is to show that this mean  converges to a distribution function which will be  described. 
Deshouillers, Dress  and Tenenbaum   [4]  studied of the analogue of this mean by considering  the set of all   divisors of  $n$, that is without  constraint on their prime factors. By denoting here by  $X_{n}$ the analogue of the random variable $X_{n,y}$ defined on the set of  all divisors of $n$, they showed that   
 \begin{eqnarray}\frac{1}{x}\sum_{n\leqslant x}\mathbb{P}(X_{n}\leqslant v) =\frac{2}{\pi}\,\arcsin(\sqrt{v})+O\left( \frac{1}{\sqrt{\log x}}\right),
 \end{eqnarray}
 uniformely for $v\in[0,1]$. This arcsin law is a Dirichlet  law in one dimension  with parameters equal $(1/2,1/2)$. Studying Couples of divisors, the authors of the present paper [7]  showed that they are distributed according  a two dimensional Dirichet Law . The method works in higher  dimension but becomes very  technique. de la Bret\`eche et Tenenbaum [2] studied  also  couples of divisors by using a probabilistic model that preserves the equiprobability of first marginal law and allows them to deduce the second marginal law. They also obtained a  Dirichlet  law.\\
Recently, Basquin [1]  studied an analogue  of this question of law of divisors by considering the set of divisors of smooth integers $n$. This question is naturally to the  de Bruijn function:
$$\Psi(x,y):= \sharp\{n\leqslant x: P(n)\leqslant y\}\cdot$$
The asymptotic behavior of de Bruijn's function is  known in a large range of  $xy$-plane. It is conneted to Dickman's function  $\rho$, which is the continuous solution in  $]0,+\infty[$  to the differential-difference  equation with inital an condition:
$$
\left
\lbrace
\begin{array}{lcl}
 u\rho^{\prime}(w)+\rho(w-1)=0,\quad (w>1)\\
 \rho(w)=1,\qquad \qquad (0\leqslant w \leqslant 1)\\
 \rho(w)=0, \quad (w<0),
\end{array}
\right.
$$
which the asymptotic behavior is well known. For example, we have
$$\log \rho (w)=-(1+0(1))w\log w, \quad (w\rightarrow +\infty).$$
Before quoting Basquin's result  and  formulate the behavior of  $\Psi(x,y)$, let us  introduce some notations that will be maintained  throughout the rest this paper. For $1< y \leqslant x$, we set
$$
u:=\frac{\log x}{\log y},
$$
and we denote by $(H_{\epsilon})$ the subset of $\mathbb{R}^{2}$ defined by the condition
$$x\geqslant x_{0}(\epsilon), \quad \exp\left((\log\log x)^{\frac{5}{3} +\epsilon}\right)\leqslant y\leqslant x,$$
where $ x_{0}(\epsilon)>0$ is a suffisantely large constant depending of  $\epsilon>0$.
Here it is  sufficient to quote the following  asymptotic formula for $\Psi(x,y)$ due to Hildebrand [6] and valid in  the range $(H_{\epsilon})$
$$\Psi(x,y)= x \rho(u)\left( 1+O\left(\frac{\log (u+1)}{\log y}\right) \right)\cdot$$
Let us introduce  functions  $\rho_{k}$  for  $k\in ]0,+\infty[$. Each function  $\rho_{k}$  is the continuous solution to the differential-difference equation with initial condition :
$$\left 
\lbrace
\begin{array}{lcl} 
w\rho'_{k}(w)+(1-k)\rho_{k}(w)+k\rho_{k}(w-1)=0, \quad  (w>1)\\
\rho_{k}(w)=\frac{1}{\Gamma(k)}\,w^{k-1},\quad (0<w \leqslant 1)\\
\rho_{k}(w)=0,\quad (w\leqslant 0).
\end{array}
\right.
$$ 
In particular, we have  $\rho_{1}=\rho$. Function  $\rho_{k}$ is the k-th fractionnal convolution  power of $\rho$ -- see Hensley's work [5].
Its asymptotic behavior  of $\rho_{k}$ is well known -- see in particular Smida's papers  [8],[9],  where we find  properties of these functions and their connection to the asymptotic behavior  of Dickman's function $\rho$. In particular, we have  the formula [8],
$$\rho_{k}(u)=k^{u(1+O(\frac{1}{\log u})}\rho(u),\quad (u\rightarrow +\infty).$$
Basquin showed that
$$\frac{1}{\Psi(x,y)}\,\sum_{n\leqslant x,P(n)\leqslant y} \mathbb{P}(X_{n}\leqslant v)=\frac{1}{\rho(u)}\,\int_{0}^{uv} \rho_{\frac{1}{2}}(s)\rho_{\frac{1}{2}}(u-s)\,ds +O\left(   \frac{\log (u+1)}{\log y}+\frac{1}{\sqrt{\log y}}       \right),$$
uniformely for  $v\in[0,1]$ and  $(x,y)\in (H_{\epsilon})$,
and he deduced  that as   $u\longrightarrow +\infty$, the distribution function converges to  the  normal distribution. More precisely, he showed that
$$\frac{1}{\rho(u)}\,\int_{0}^{uv} \rho_{\frac{1}{2}}(s)\rho_{\frac{1}{2}}(u-s)\,ds=\Phi\left(u\sqrt{2\xi'(u)}(v-\frac{1}{2})\right)+O\left(\frac{1}{u}\right),$$
where  $\xi'(u)\sim 1/u $, as $u\rightarrow +\infty$ and 
$$\Phi(w)=\frac{1}{\sqrt{\pi}} \int_{-\infty}^{w} e^{-t^{2}}\,dt,$$
is the normal distribution function.
\section{Statements of results} 
To state our results, let us introduce the Buchstab function  $\omega$. This function was discovered by Bushstab  and comes from the study of the uncancelled elements in the sieve of Eratosthenes -- see de Bruijn's  beautiful paper [3]. It is the unique continuous solution for  $v>1$ to the  differential-difference equation, with initial condition:
$$
\left \lbrace
\begin{array}{lcl}
v\omega'(v)+\omega(v) -\omega(v-1)=0,\quad (v>2)\\
\omega(v)=\frac{1}{v},\quad (1\leqslant v\leqslant 2)\\
\omega(v)=0,\quad v<1.
\end{array}
\right.
$$
Its asymptotic behavior is known -- see de Bruijn [3] and Tenenbaum's book  [10 ] chap. III.6). In particular for $v\geqslant0$, we have
\begin{eqnarray}
\omega(v)=e^{-\gamma}+O\left(\rho(v)e^{\frac{-c\,v}{\log^{2}(v+2)}}\right),
\end{eqnarray}
where  $\gamma$ is the Euler function and $c$ a positive constant.
In the first theorem below we show the convergence  of the mean of distribution functions  (1) to a distribution function,
in the second one we describe the limit law as $u\longrightarrow  +\infty$ and in the third one, we give as an example  expressions of  the limit law for
$1\leqslant u \leqslant  2$. 
\begin{Theo}
Uniformly for $v\in[0,1]$ and 
 $(x, y)$ in $(H_{\epsilon})$, we have
 $$
\begin{array}{lcl}
\frac{1}{x}\sum_{n\leqslant x}\mathbb{P}(X_{n,y}\leqslant v)&=&\int_{0}^{uv} \left(\int_{0}^{u-s-1}\rho_{\frac{1}{2}}(z)\omega(u-s-z )\,dz\right)\rho_{\frac{1}{2}}(s)\,ds \\[0.5cm]
&+&\int_{0}^{uv} \rho_{\frac{1}{2}}(s)\,\rho_{\frac{1}{2}}(u-s)\,ds\\[0.5cm]
&+&O\left( \frac{\log (u+1)}{\log y}+\frac{1}{\sqrt{\log y}}\right)\cdot
\end{array}
$$
\end{Theo}
We notice that for $y=x$, that is to say $u=1$, the  formula of Theorem 2.1  is reduced to formula (2) obtained  in [4]. Indeed, 
the first integral vanishes  because  $\rho_{\frac{1}{2}}(z)=0$ for $z\leqslant 0$ and
$$\int_{0}^{v} \rho_{\frac{1}{2}}(s)\,\rho_{\frac{1}{2}}(1-s)\,ds=\frac{2}{\pi}\int_{0}^{\sqrt{v}}\,\frac{dt}{\sqrt{1-t^{2}}}=\frac{2}{\pi}\arcsin(\sqrt{v}).$$
Let us denote 
$$F(u,v)=\int_{0}^{uv} \left(\int_{0}^{u-s-1}\rho_{\frac{1}{2}}(z)\omega(u-s-z )\,dz\right)\rho_{\frac{1}{2}}(s)\,ds+\int_{0}^{uv} \rho_{\frac{1}{2}}(s)\,\rho_{\frac{1}{2}}(u-s)\,ds.$$
We have 
\begin{Theo}
 For  $v\in[0,1]$ and  as $u\rightarrow +\infty$, we uniformly have
 $$F(u,v)=\frac{1}{\sqrt{e^{\gamma}}}  \int_{0}^{uv}\rho_{\frac{1}{2}}(s)ds+O\left(\rho_{2}(u) \right),$$
 where $\gamma$ is the  d'Euler constant.
\end{Theo}
\begin{Theo}
\begin{enumerate}
\item
For $v\in[0,\frac{u-1}{u}]$ and $1<u\leqslant 2$, we have 
$$
 F(u,v)=\frac{2}{\pi}\arcsin(\sqrt{v})+\frac{1}{\pi}\left(\;\log(u)+\log(1-v)\;\right)\arcsin\left(  \sqrt{\frac{uv}{u-1}} \,  \right)-\frac{1}{2}\log(1-v).
 $$
 \item
 For $v\in[\frac{u-1}{u},\frac{1}{u}]$ and $1<u\leqslant 2$, we have 
 $$F(u,v)=\frac{2}{\pi}\arcsin(\sqrt{v})+\frac{1}{2}\log u.$$
 \item
 For $v\in[\frac{1}{u},1]$ and $1<u\leqslant 2$, we have
 $$F(u,v)=\frac{2}{\pi}\arcsin(\sqrt{v})+\frac{1}{\pi}\left(\;\log(u)+\log(v)\;\right)\arcsin\left(  \sqrt{\frac{u(1-v)}{u-1}} \,  \right)-\frac{1}{2}\log(v)\cdot$$
 \end{enumerate}
\end{Theo}
Let us set
$$\begin{array}{lcl}
S(x,y,v)&:=&\sum_{n\leqslant x}\mathbb{P}(X_{n,y}\leqslant v)=\sum_{n\leqslant x}\frac{1}{\tau(n,y)}\,\sum_{d|n, d\leqslant n^{v},P(d)\leqslant y} 1\\[0.5cm]
&=& S_{1}(x,y,v)-S_{2}(x,y,v),
\end{array}
$$
with
$$S_{1}(x,y,v):= \sum_{n\leqslant x}\frac{1}{\tau(n,y)}\,\sum_{\stackrel{d|n, d\leqslant x^{v}}{P(d)\leqslant y}}1=\sum_{d\leqslant x^{v}, P(d)\leqslant y}\,\sum_{m\leqslant x/d}\frac{1}{\tau(dm,y)}, $$
and 
$$S_{2}(x,y,v):= \sum_{n\leqslant x}\frac{1}{\tau(n,y)}\,\sum_{\stackrel{d|n, n^{v}<d\leqslant x^{v}}{P(d)\leqslant y}}1.$$
We will show that the main contribution to Theorem 2.1 comes from the estimation of  $S_{1}(x,y,v)$. The proof rests on the estimation 
$$\sum_{n\leqslant x}\frac{1}{\tau(dn,y)}\; ,$$
for $(x,y)$ in $(H_{\epsilon})$ and $d\geqslant 1$, $y$-smooth. To study  this quantity, we used the factorisation  $n=ab$ with  $a$ is $y$-smooth and  $b$  has all its prime factors greater than $y$. This allowed us to use results on this topic  available in the literature. 
\section{Proof of  Theorem 2.1 } 
\subsection{Preparatory Lemmas}
Let us introduce some notations which will be used in the sequel. For each fixed integer $d\geqslant 1$ we define a multiplicative function 
$$\gamma_{d}(n)=\frac{\tau(d)}{\tau(dn)},$$
where $\tau$ is the divisor function.
For each prime number $p$, we denote  by $v_{p}(d)$ the p-adic valuation of $d$. We have
 $$\gamma_{d}(p^{\alpha})=\frac{v_{p}(d)+1}{v_{p}(d)+\alpha+1}\cdot$$
 we consider  the Dirichlet series of $\gamma_{d}(n)$
$$F_{d}(s):=\sum_{n\geqslant1}\frac{\gamma_{d}(n)}{n^{s}}, \quad ( \Re (s)>1).$$ We have
  $F_{d}(s)=\zeta^{\frac{1}{2}}(s)G_{d}(s)$ the half-plane  $\Re (s)>1$, 
with $\zeta$ is the Riemann zeta function and
 $$
 \begin{array}{lcl}
 G_{d}(s)&=&\prod_{p}\left( 1-\frac{1}{p^{s}}\right)^{\frac{1}{2}}\left( \sum_{\alpha\geqslant 0}\frac{v_{p}(d)+1}{(v_{p}(d)+\alpha+1)p^{\alpha s}}\right) \\[0.5cm]
 &=& B(s) K_{d}(s),
 \end{array}
 $$
 with 
  $$B(s):=\prod_{p}\left(1-\frac{1}{p^{s}}\right)^{\frac{1}{2}}\left(1+\sum_{\alpha \geqslant 1}\frac{1}{(\alpha+1)p^{\alpha s}}\right),$$ 
 and
 $$K_{d}(s):=\prod_{p^{\beta}\| d}\left(1+ \sum_{\alpha\geqslant 1}\frac{\beta+1}{(\beta+\alpha+1)p^{\alpha s}}\right)\left(1+\sum_{\alpha\geqslant 1}\frac{1}{(\alpha+1)p^{\alpha s}}\right)^{-1}\cdot$$
 For each fixed integer $d\geqslant1$, we define Dirichlet series
$$K_{d}(s):=\sum_{n\geqslant 1}\frac{\delta_{d}(n)}{n^{s}}; \quad B(s)=\sum_{n\geqslant 1}\frac{b(n)}{n^{s}};\quad G_{d}(s)=\sum_{n\geqslant 1}\frac{h_{d}(n)}{n^{s}}\cdot$$
Then we have
\begin{eqnarray}
F_{d}(s)=\sum_{n\geqslant1}\frac{\gamma_{d}(n)}{n^{s}}=\zeta^{\frac{1}{2}}(s)\sum_{n\geqslant 1}\frac{h_{d}(n)}{n^{s}}
\end{eqnarray}
and 
\begin{eqnarray}
\sum_{n\geqslant 1}\frac{h_{d}(n)}{n^{s}}=\left( \sum_{n\geqslant 1}\frac{b(n)}{n^{s}}\right)\left(\sum_{n\geqslant 1}\frac{\delta_{d}(n)}{n^{s}}\right)\cdot
\end{eqnarray}
\begin{Lem}
Let $d\geqslant1$ be a fixed integer. It exists a real number  $0<\eta< 1/3$ such that the series
$$\sum_{n\geqslant 1}\frac{h_{d}(n)}{n^{s}}$$
is absolutely  convergent in the half-plane $\Re(s)=\sigma\geqslant1-\eta$ and for $d\geqslant1$  we uniformly have  
$$\sum_{n\geqslant 1}\frac{h_{d}(n)}{n^{\sigma}}\ll_{\eta} \prod_{p|d}(1+\frac{2}{p^{\sigma}})\cdot$$
\end{Lem}
The  result of this lemma can be deduced from from a general study developped in dans [4],  page 275. Indeed, we apply  lemme 1 of [4], page 276, we get that the series 
$$\sum_{n\geqslant 1}\frac{\delta_{d}(n)}{n^{s}}=K_{d}(s)$$is absolutely convergent in the half-plane   $\Re(s)=\sigma\geqslant1-\eta$  and we have 
\begin{eqnarray}\sum_{n\geqslant 1}\frac{|\delta_{d}(n)|}{n^{\sigma}}\ll_{\eta} \prod_{p|d}(1+\frac{2}{p^{\sigma}})\cdot
\end{eqnarray}
The  lemme 2 of  [4], page 278, applies to the series 
$$\sum_{n\geqslant 1}\frac{b(n)}{n^{s}}=B(s)$$
with  the exponent  $\alpha=1/2$ and $\psi(n)=1/\tau(n)$. We obain 
 \begin{eqnarray}\sum_{n\geqslant 1}\frac{|b(n)|}{n^{\sigma}}\leqslant \prod_{p}( 1+ O(\frac{1}{p^{2-2\eta}}))\ll_{\eta}1.
 \end{eqnarray}
 in the half-plane $\Re(s)=\sigma\geqslant1-\eta$. Lemma 3.1 follows from  $(5)$, $(6)$ and $(7)$.
 \hfill{
$\square$}\\[0.5cm]
We set
 $$M_{\eta,d}:=\prod_{p|d}(1+\frac{2}{p^{1-\eta}})\cdot$$
From the formula $(4)$ we define a mutiplicative function  $h_{d}$ by the convolution identity $\gamma_{d}=\tau_{\frac{1}{2}}\ast h_{d}$. 
From lemma 3.1 we get that its Dirichlet  series satisfies the conditions (1.18) of  Th\'eor\`eme 3 of  [9], page 25,  since we have  
$$\sum_{\stackrel{n>t}{P(d)\leqslant y}}\frac{|h_{d}(n)|}{n}\leqslant\frac{1}{t^{\eta/2}}\sum_{\stackrel{n>t}{P(d)\leqslant y}}\frac{|h_{d}(n)|}{n^{1-(\eta)/2})}\leqslant\frac{1}{t^{\eta/2}}\sum_{n\geqslant1}\frac{|h_{d}(n)|}{n^{1-(\eta/2)}}\ll_{\eta} \frac{M_{\eta,d}}{t^{\eta/2}}\cdot$$
The proof of Théorème 3 of [9]  paragraph 6 page 45  works  and we obtain the first result of  lemma 3.2 below. The second result is a consequence of   Théorème T of [4]  applied to the function  $1/\tau(n)$, and by noting that for $0<u\leqslant 1$  
 \begin{eqnarray}
 \rho_{\frac{1}{2}}(u)=\frac{1}{\sqrt{\pi} \sqrt{u}}= \frac{\sqrt{\log y}}{\sqrt{\pi}\,\sqrt{\log x}}\cdot
 \end{eqnarray}
  \begin{Lem} 
 \begin{enumerate}
 \item
  let be $\eta\in]0,\frac{1}{2}[$ and let  $\epsilon>0$  be fixed . For $d\geqslant1$ and $(x,y) \in (H_{\epsilon})$ we uniformly have
$$\sum_{\stackrel{n\leqslant x}{P(n)\leqslant y}}\gamma_{d}(n)=\frac{x}{\sqrt{\log y} }\;\rho_{\frac{1}{2}}(u)\left( G_{d}(1)+O\left(M_{\eta,d}\,\left(\frac{\log (u+1)}{\log y}+\frac{1}{\sqrt{ \log y}}\right)\right) \right)\cdot$$

\item
For  $1<x\leqslant y$ and $d\geqslant1$, we uniformly have
$$\sum_{\stackrel{n\leqslant x}{P(n)\leqslant y}}\gamma_{d}(n)=\frac{x}{\sqrt{\log y} }\;\rho_{\frac{1}{2}}(u)\left( G_{d}(1)+O\left(\frac{M_{\eta,d}}{\log x}\right) \right)\cdot$$
\end{enumerate}
\end{Lem}
  Basquin [2] obtained the first result of this lemma by using another  convolution identity and by applying a general result of  Tenenbaum et Wu [11]. We set
$$g(d):=\frac{K_{d}(1)}{\tau(d)}=\prod_{p^{\beta}\| d}\left(\sum_{\alpha\geqslant 0}\frac{1}{(\beta+\alpha+1)p^{\alpha }}\right)\left(\sum_{\alpha\geqslant 0}\frac{1}{(\alpha+1)p^{\alpha }}\right)^{-1}\cdot$$
$g$ is a multiplicative function and for  $\Re(s)>1$, we have 
$$\sum_{n\geqslant 1}\frac{g(n)}{n^{s}}=\zeta^{\frac{1}{2}}(s)\,\sum_{n\geqslant 1}\frac{\beta(n)}{n^{s}},$$
where $\beta$ is a multiplicative function satisfying $g=\tau_{\frac{1}{2}}\ast \beta$. We have 
$$H(s):=\sum_{n\geqslant 1}\frac{\beta(n)}{n^{s}}=\prod_{p}\left(1-\frac{1}{p^{s}}\right)^{\frac{1}{2}}\left( \sum_{\alpha\geqslant0} \frac{g(p^{\alpha})}{p^{\alpha s}} \right)\cdot$$
Let us note for later use that
\begin{eqnarray}
G_{d}(1)=K_{d}(1)B(1)=\tau(d) g(d) B(1).
\end{eqnarray}
 Lemme 2 of  [4], page 278, applies to the series  $H(s)$
 with  exponent  $\alpha=1/2$ and the function  $\psi(n)=g(n)$. It follows that the series $H(s)$  is absolutely convegent  in the half-plane  $\Re(s)=\sigma\geqslant1-\eta$ and we have
$$\sum_{n\geqslant 1}\frac{|\beta(n)|}{n^{\sigma}} \ll_{\eta}1.$$
 The conditions (1.18) of application of 
Th\'eor\`eme 3 of  [9], page 25, are satisfied and we obtain the first result of the following lemma 3.3. The second result is an immediate consequence of lemme 3 page 282 of [4] and  relation (8) above.
\begin{Lem} 
\begin{enumerate}
\item
 Uniformly in the range $(H_{\epsilon})$. We have
$$\sum_{\stackrel{n\leqslant x}{P(n)\leqslant y}}g(n)=H(1)\,\frac{x}{\sqrt{\log y} }\;\rho_{\frac{1}{2}}(u)\left( 1+O\left(\frac{\log (u+1)}{\log y}+\frac{1}{\sqrt{ \log y}}\right) \right)\cdot$$
 
\item
For $1<x\leqslant y$, we have 
$$\sum_{n\leqslant x}g(n)=H(1)\,\frac{x}{\sqrt{\log y} }\;\rho_{\frac{1}{2}}(u)\left( 1+O\left(\frac{1}{\log x } \right)\right)\cdot$$
\end{enumerate}
\end{Lem}
 \begin{Rem}
 We set 
$$N(d):=\frac{M_{\eta,d}}{\tau(d)}\cdot$$
The function $N$ is multiplicative and positive and satisfies
$N(p)=1+O(\frac{1}{\sqrt{p}})$ and $N(p^{\alpha})\ll 1$. 
So,  from Lemme 3.1 of  Basquin [1],  and partial summation we get
$$\sum_{\stackrel{d\leqslant x}{P(d)\leqslant y}}\frac{N(d)}{d}\ll_{\eta} \sqrt{\log y}\,,$$
uniformly in  $(H_{\epsilon})$.
\end{Rem}
\begin{Lem}
We have $B(1)H(1)=1$.
\end{Lem}
\textbf{Proof.}
We have 
$$B(1)=\prod_{p}\left(1-\frac{1}{p}\right)^{\frac{1}{2}}\left(\sum_{j \geqslant 0}\frac{1}{(j+1)p^{j}}\right),\quad H(1)=\prod_{p}\left(1-\frac{1}{p}\right)^{\frac{1}{2}}\left(\sum_{\alpha\geqslant0}\frac{g(p^{\alpha})}{p^{\alpha}}\right),$$
and
$$ \sum_{\alpha\geqslant0}\frac{g(p^{\alpha})}{p^{\alpha}}= \sum_{\alpha\geqslant0}\left(  \sum_{j \geqslant 0}\frac{1}{(\alpha +j+1)p^{\alpha+j}} \right)\left(  \sum_{j \geqslant 0}\frac{1}{(j+1)p^{j}}  \right)^{-1},$$
then
$$B(1)H(1)=\prod_{p}\left(1-\frac{1}{p}\right)\sum_{\alpha\geqslant0}\left(  \sum_{j \geqslant 0}\frac{1}{(\alpha +j+1)p^{j+\alpha}} \right)\cdot$$
By noticing that for $0<|x|<1$,
$$\frac{x^{\alpha +j}}{\alpha +j+1 }=\frac{1}{x}\int_{0}^{x}t^{\alpha +j} dt,$$
 we get
$$\sum_{\alpha\geqslant0}\left(  \sum_{j \geqslant 0}\frac{x^{\alpha +j}}{\alpha +j+1} \right)=  \frac{1}{x}\int_{0}^{x}\frac{dt}{(1-t)^{2}}=\frac{1}{1-x}, $$
and then  $B(1)H(1)=1$.
\hfill{
$\square$}\\[0.5cm]
Let us  denote by $P_{-}(n)$ the smallest prime factors of $n>1$ and we set  $P_{-}(1)=+\infty$. 
 The function
$$ \Phi(x,y):= \sum_{\stackrel{n\leqslant x}{P_{-}(n)>y}}1$$
 has been studied by  de Bruijn [3]. In the lemma below we quote an asymptotic formula which is sufficient for our purpose. 

\begin{Lem}
Uniformly in the range $(H_{\epsilon})$ we have 
$$\Phi(x,y)=\frac{x\,\omega(u)}{\log y} -\frac{y}{\log y}+O\left(\frac{x}{(\log y)^{2}}\right).$$

\end{Lem}
\textbf{Proof.} Lemma follows from  Corollaire 6.14 of [10] chap. III.6 page 547 and Mertens Formula. \hfill{
$\square$}
\begin{Lem} 
 Uniformly for each integer $d$ such that  $P(d)\leqslant y $ and $ (x,y)$  in $(H_{\epsilon})$ we have
$$
\begin{array}{lcl}
\sum_{n\leqslant x}\frac{1}{\tau(dn,y)}&=&\frac{x}{\sqrt{\log y}}\,\frac{G_{d}(1)}{\tau(d)}\left(\int_{0}^{u-1}\rho_{\frac{1}{2}}(z)\omega(u-z)\,dz +\rho_{\frac{1}{2}}(u)\right)\\[0.5cm]
&+& O\left( \frac{M_{\eta,d}}{\tau(d)} \left(\frac{x\,\log (u_1)}{(\log y)^{\frac{3}{2}}}+\frac{x}{\log y}\right)\right)\cdot$$
\end{array}$$
\end{Lem}
\textbf{Proof.} We write $n=ab$ with $P(a)\leqslant y $ and ($b=1$ or  if $p| b\Rightarrow p>y $).   Then for $d$ $y$-smooth, we have
$$\tau(dn,y)=\tau(dab,y)=\tau(da).$$
and then
$$T_{d}(x,y):=\sum_{n\leqslant x}\frac{1}{\tau(dn,y)} =\sum_{a\leqslant x, P(a)\leqslant y}\frac{1}{\tau(da)}\, \Phi(\frac{x}{a},y),$$
First we consider the range $$\exp\left((\log\log x)^{\frac{5}{3} +\epsilon}\right)\leqslant y\leqslant \frac{x}{a}\; ,\quad  x\geqslant x_{0}(\epsilon).$$ 
Write
$$\sum_{n\leqslant x}\frac{1}{\tau(dn,y)} = \sum_{a\leqslant x/y, P(a)\leqslant y}\frac{1}{\tau(da)}\, \Phi(\frac{x}{a},y)+\sum_{x/y<a\leqslant x, P(a)\leqslant y}\frac{1}{\tau(da)} :=T_{d}(x,y)+\overline{T}_{d}(x,y).$$
Lemma 3.2 gives
$$\overline{T}_{d}(x,y)=\frac{G_{1}(d)}{\tau(d)}\;\frac{x}{\sqrt{\log y}}\,\rho_{\frac{1}{2}}(u)+O\left( \frac{M_{\eta,d}}{\tau(d)} \left(\frac{x\,\log (u+1)}{(\log y)^{\frac{3}{2}}}+\frac{x}{\log y}\right)\right)\cdot$$
To study  $T_{d}(x,y)$ we apply Lemma 3.5. We obtain
$$
\begin{array}{lcl}
T_{d}(x,y)&=&\frac{x}{\log y}\sum_{a\leqslant x/y, P(a)\leqslant y}\frac{1}{a\tau(da)}\omega(u-\frac{\log a}{\log y})+ O\left( \frac{x}{\log^{2} y}\sum_{a\leqslant x/y, P(a)\leqslant y}\frac{1}{a\tau(da)} \right)\\
&=:& \frac{x}{\log y}\; L_{1}+O\left( \frac{x}{\log^{2} y} \;L_{2}   \right)\cdot
\end{array}$$
Now let us study $L_{1}$.  Partial summation and Lemma 3.2 give
$$\begin{array}{lcl}
L_{1}
&=&\frac{1}{\tau(d)}\int_{1^{-}}^{x/y}\frac{1}{t}\omega(u-\frac{\log t}{\log y})\,d\left(\sum_{a\leqslant t, P(a)\leqslant y}\gamma_{d}(a)\right)\\[0.5cm]
&=&\frac{1}{\tau(d)}\int_{1}^{x/y}\left(\sum_{a\leqslant t, P(a)\leqslant y}\gamma_{d}(a)\right) \,\left(\frac{1}{t^{2}}\omega(u-\frac{\log t}{\log y})\,+\omega'(u-\frac{\log t}{\log y})\frac{1}{t^{2} \log y}\right)\,dt\\[0.5cm]
&+&O\left( \frac{M_{\eta,d}}{\tau(d)} \frac{1}{\sqrt{\log y}}\right)\\[0.5cm]
&=&\frac{1}{\sqrt{\log y}}\frac{G_{d}(1)}{\tau(d)}\int_{1}^{x/y}\left(\frac{1}{t}\omega(u-\frac{\log t}{\log y})\,+\omega'(u-\frac{\log t}{\log y})\frac{1}{t \log y}\right)\,\rho_{\frac{1}{2}}(\frac{\log t}{\log y})\,dt\\[0.5cm]
&+& O\left(\frac{M_{\eta,d}}{\tau(d)} \int_{1}^{x/y}\left(\frac{\log (\frac{\log t}{\log y}+1)}{(\log y)^{\frac{3}{2}}}+\frac{1}{ \log y} \right) \omega(u-\frac{\log t}{\log y})\rho_{\frac{1}{2}}(\frac{\log t}{\log y})\frac{dt}{t}\right)\\[0.5cm]
&+&O\left(\frac{M_{\eta,d}}{\tau(d)} \int_{1}^{x/y}\left(\frac{\log (\frac{\log t}{\log y}+1)}{(\log y)^{\frac{3}{2}}}+\frac{1}{ \log y} \right) |\omega'(u-\frac{\log t}{\log y})|\rho_{\frac{1}{2}}(\frac{\log t}{\log y})\frac{dt}{t \log y}\right)\\[0.5cm]
&+&O\left( \frac{M_{\eta,d}}{\tau(d)} \frac{1}{\sqrt{\log y}}\right)\cdot 
\end{array} $$
 By  change of variable $z=\frac{\log t}{\log y}$. We obtain
$$
\begin{array}{lcl}
L_{1}
 &=&\sqrt{\log y}\,\frac{G_{d}(1)}{\tau(d)}\int_{0}^{u-1}\rho_{\frac{1}{2}}(z)\omega(u-z)\,dz + O\left( \frac{M_{\eta,d}}{\tau(d)} \left(1+\frac{\log u}{(\log y)^{\frac{1}{2}}}\right)\right),
\end{array}
$$
 We estimate $L_{2}$ in the same way. Lemma 3.2 and By partial summation give
$$L_{2}\ll\sqrt{\log y}\;\frac{M_{\eta,d}}{\tau(d)} \int_{0}^{u}\rho_{\frac{1}{2}}(z)dz\ll\sqrt{\log y}\;\frac{M_{\eta,d}}{\tau(d)}\cdot$$
We get our result in the considered range  from these different estimates. In the range  $\frac{x}{a} <y\leqslant x$  we have 
$\Phi(x/a,y)=1$. The result follows immediately from Lemma 3.2.
\hfill{
$\square$}
\begin{Lem} 
  Set  $\epsilon_{x}=\frac{\log 2 }{\log x}$. Uniformly for $(x,y)\in (H_{\epsilon})$ \\and $0\leqslant v < \epsilon_{x}$ we have
 $$ S(x,y,v)=O\left(\frac{x}{\sqrt{\log y}}\right)\cdot$$
 \end{Lem}
\textbf{Proof} The condition on  $v$ implies  $d=1$. By   Lemme 3.6 with  $d=1$ we get
$$S(x,y,v)=\sum_{n\leqslant x}\frac{1}{\tau(n,y)}\,\sum_{d|n, d\leqslant n^{v},P(d)\leqslant y} 1 \leqslant \sum_{n\leqslant x}\frac{1}{\tau(n,y)}\ll \frac{x}{\sqrt{\log y }}\cdot$$
\hfill{
$\square$}
\begin{Lem} Let be $\epsilon_{x}=\frac{\log 2 }{\log x}$. Uniformly for  $(x,y)\in (H_{\epsilon})$ and $ \epsilon_{x}\leqslant v\leqslant1$, we have
$$S_{2}(x,y,v)=O\left(\frac{x}{\log y }  \right)\cdot$$
\end{Lem}
\textbf{Proof} We have
$$S_{2}(x,y,v)=\sum_{n\leqslant x}\frac{1}{\tau(n,y)}\,\sum_{d|n, n^{v}<d\leqslant x^{v},P(d)\leqslant y} 1\leqslant \sum_{d\leqslant x^{v},P(d)\leqslant y}\,\sum_{m\leqslant x^{1-v}}\frac{1}{\tau(dm,y)}\cdot$$

If $1-\epsilon_{x}\leqslant v\leqslant1$, then $x^{1-v}<2$. In this case, we have

$$S_{2}(x,y,v)\leqslant \sum_{d\leqslant x,P(d)\leqslant y}\,\frac{1}{\tau(d)}\ll\frac{x}{\sqrt{\log y}}\,\rho_{\frac{1}{2}}(u)\ll\frac{x}{\sqrt{\log y}}\cdot$$
Suppose that  $\epsilon_{x}\leqslant v\leqslant 1-\epsilon_{x}$.
If   $y\leqslant \min\{ x^{v},x^{1-v}\}$.  We apply Lemmas 3.6 and 3.3,  we get
$$
\begin{array}{lcl}
S_{2}(x,y,v)&\ll& \sum_{d\leqslant x^{v},P(d)\leqslant y}\frac{x^{1-v}}{\sqrt{\log y}}\,\frac{G_{d}(1)}{\tau(d)}\left(\int_{0}^{(1-v)u-1}\rho_{\frac{1}{2}}(z)\omega((1-v)u-z)\,dz+\rho_{\frac{1}{2}}((1-v)u)\right) \\[0.8cm]
 &\ll& \frac{x^{1-v}}{\sqrt{\log y}}\,\sum_{d\leqslant x^{v},P(d)\leqslant y} g(d)\\ [0.8cm]
 &\ll&\frac{x}{\log y} \,\rho_{\frac{1}{2}}(uv)\ll \frac{x}{\log y},
\end{array}
$$
since $uv\geqslant1$ et $(1-v)u\geqslant1$.
 We proceed in the same way in other cases. For  $x^{1-v}< y \leqslant x^{v}$, the inner sum is
$$\sum_{m\leqslant x^{1-v}}\frac{1}{\tau(dm,y)}=\sum_{m\leqslant x^{1-v}}\frac{1}{\tau(dm)}\ll \frac{x^{1-v}}{\sqrt{(1-v)\log x}} \ll x^{1-v}$$ 
from Lemma  3.2.  We apply Lemme 3.3 to the outer sum we get
$$S_{2}(x,y,v)\ll \frac{x}{\sqrt{\log y}}\rho_{\frac{1}{2}}(uv)\ll \frac{x}{\sqrt{\log y}}\cdot$$
The case $x^{v}< y \leqslant x^{1-v}$ is similar. Last  if $y>\max\{x^{1-v},x^{v}\}$, we set
$\epsilon'_{x}=\frac{\log 2}{\sqrt{\log x}}$ and consider  three cases  $\epsilon_{x}\leqslant v\leqslant \epsilon'_{x}$ , $ \epsilon'_{x}< v\leqslant 1-\epsilon'_{x}$ et $1-\epsilon'_{x}< v \leqslant 1-\epsilon_{x}$.
In each situation we have $\frac{1}{\sqrt{v(1-v)}}\ll\sqrt{\log x}$.
 By applying Lemmas 3.2 and 3.3 we get
$$ S_{2}(x,y,v)\ll \frac{x}{\log x} \frac{1}{\sqrt{v(1-v)}}\ll \frac{x}{\sqrt{\log x}}\cdot$$
 \hfill{
$\square$}
\subsection{Proof of Theorem 2.1 } 
Taking into account Lemmas 3.7 and  3.8,  it remains to estimate $S_{1}(x,y,v)$ for  
$\epsilon_{x}\leqslant v\leqslant1$.
We consider two following situations : $\epsilon_{x}\leqslant v\leqslant1-\epsilon_{x}$ et $1-\epsilon_{x}<v\leqslant1$.
 First suppose that $1-\epsilon_{x}< v\leqslant1$. Dans ce cas $x/2 <x^{v}$. We write
 $$
 \begin{array}{lcl}
 S_{1}(x,y,v)&=&\sum_{d\leqslant x/2, P(d)\leqslant y}\,\sum_{m\leqslant x/d}\frac{1}{\tau(dm,y)}+ \sum_{x/2<d\leqslant x^{v}, P(d)\leqslant y}\,\sum_{m\leqslant x/d}\frac{1}{\tau(dm,y)}\\[0.8cm]
 &=&:\widehat{S}_{1} +\widehat{S}_{2}.
 \end{array}$$
 Let us study $ \widehat{S}_{2}$. $ x/2< d \iff x/d<2$ thus $m=1$. We have
 $$ \widehat{S}_{2}=\sum_{x/2<d\leqslant x^{v}, P(d)\leqslant y}\frac{1}{\tau(d)}\ll\sum_{d\leqslant x, P(d)\leqslant y}\frac{1}{\tau(d)}\ll\frac{x}{\sqrt{\log y}}\rho_{\frac{1}{2}}(u)\ll\frac{x}{\sqrt{\log y}}\cdot$$
The evalation of $\widehat{S}_{1}$ is similar to the envaluation of  $S_{1}(x,y,v)$ under  complementary condition $\epsilon_{x}\leqslant v\leqslant 1-\epsilon_{x}$, sinc  $ x/2=x^{1-\epsilon_{x}}$.
Let us study  $S_{1}(x,y,v)$ under the condition $\epsilon_{x}\leqslant v\leqslant 1-\epsilon_{x}$. First, consider the range $$\exp\left((\log\log x)^{\frac{5}{3} +\epsilon}\right)\leqslant y\leqslant x/d\; ,\quad  x\geqslant x_{0}(\epsilon).$$
We write
 $$
 \begin{array}{lcl}
 S_{1}(x,y,v)&:=&\sum_{n\leqslant x}\frac{1}{\tau(n,y)}\,\sum_{\stackrel{d|n, d\leqslant x^{v}}{P(d)\leqslant y}}1=\sum_{d\leqslant x^{v}, P(d)\leqslant y}\,\sum_{m\leqslant x/d}\frac{1}{\tau(dm,y)}\\[0.5cm]
&=&\sum_{d\leqslant x^{v}, P(d)\leqslant y}\,T_{d}(\frac{x}{d},y)\cdot
\end{array}
$$
We apply Lemma 3.6.  We obtain
$S_{1}(x,y,v)=S_{1,1}+O(S_{1,2})$ avec, en posant  $$u_{d}:=u-\frac{\log d}{\log y},$$
$$S_{1,1}:=\frac{x}{\sqrt{\log y}}\,\sum_{d\leqslant x^{v}, P(d)\leqslant y}\frac{G_{d}(1)}{d\tau(d)}\left(\int_{0}^{u_{d}-1}\rho_{\frac{1}{2}}(z)\omega(u_{d}-z)\,dz+\rho_{\frac{1}{2}}(u_{d})\,\right),$$
 and
$$S_{1,2}:=x\,\left(\frac{1}{\log y}+\frac{\log(u+1)}{(\log y)^{\frac{3}{2}}}\right)\,\sum_{d\leqslant x^{v}, P(d)\leqslant y}\frac{M_{\eta,d}}{d\tau(d)}\ll x\,\left(\frac{1}{\sqrt{\log y}}+\frac{\log(u+1)}{\log y}\right), $$
by using Remark 3.1.
It remains to estimate  $S_{1,1}$. From (9), we have  $\frac{G_{d}(1)}{\tau(d)}=B(1) g(d)$.
 By partial summation we get
 $$
\begin{array}{lcl} 
 S_{1,1}&=&\frac{x}{\sqrt{\log y}}\,B(1)\sum_{ d\leqslant x^{v}, P(d)\leqslant y}\frac{g(d)}{d}\left(\int_{0}^{u_{d}-1}\rho_{\frac{1}{2}}(z)\omega(u_{d}-z)\,dz+\rho_{\frac{1}{2}}(u_{d})\right)\\
 &=&\frac{B(1)x^{1-v}}{\sqrt{\log y}}\left(\int_{0}^{(1-v)u-1} \rho_{\frac{1}{2}}(z)\omega((1-v)u-z )\, dz+\rho_{\frac{1}{2}}((1-v)u) \right)\left(\sum_{n\leqslant x^{v}, P(n)\leqslant y}g(n)\right)\\
 &-&\frac{x}{\sqrt{\log y}}\,B(1) \;\int_{1}^{x^{v}}
\left( \sum_{n\leqslant t, P(n)\leqslant y}g(n)\right)\;d\left(\frac{1}{t}\left(\int_{0}^{u_{t}-1} \rho_{\frac{1}{2}}(z)\omega(u_{t}-z )\, dz+\rho_{\frac{1}{2}}(u_{t})\right)\right)\\
&:=& R_{1}-R_{2}.
 \end{array}
 $$
From  Lemmas 3.3 and  3.4  
 $$
R_{1}\ll\frac{x}{\log y}\, \rho_{\frac{1}{2}} (uv)\left(1+ \rho_{\frac{1}{2}}((1-v)u)\right)
\ll \frac{x}{\sqrt{\log y}}\cdot
$$
the last upper bound follows by using the same way as in proof of Lemma 3.8.
Now consider  $R_{2}$.  We have
 $$d\left( \frac{1}{t}  \left(\int_{0}^{u_{t}-1}\rho_{\frac{1}{2}}(z)\,\omega(u_{t}-z )\, dz+\rho_{\frac{1}{2}}(u_{t}) \right)\right) =D_{1}+D_{2},$$
 with 
 $$D_{1}:=-\left(\int_{0}^{u_{t}-1}\,\rho_{\frac{1}{2}}(z)\omega(u_{t}-z )\,dz +\rho_{\frac{1}{2}}(u_{t})\right)\,\frac{dt}{t^{2}},$$
and
$$D_{2}:=\frac{1}{t}d\left( \int_{0}^{u_{t}-1}\rho_{\frac{1}{2}}(z)\,\omega(u_{t}-z )\, dz +\rho_{\frac{1}{2}}(u_{t})\right).$$
 Write $ R_{2}=R_{2,1}+R_{2,2}$
with
$$R_{2,1}:=\frac{x}{\sqrt{\log y}}\,B(1) \;\int_{1}^{x^{v}}  \left(\sum_{n\leqslant t, P(n)\leqslant y}g(n)\right) D_{1},$$
and
$$R_{2,2}:=\frac{x}{\sqrt{\log y}}\,B(1) \;\int_{1}^{x^{v}}  \left(\sum_{n\leqslant t, P(n)\leqslant y}g(n)\right) D_{2}.$$
To evaluate  $R_{2,1}$ we apply Lemmas 3.3 and  3.4. We get
$$
\begin{array}{lcl}
R_{2,1}&=&
-\frac{x}{\log y}\int_{1}^{x^{v}}\rho_{\frac{1}{2}}(\frac{\log t}{\log y})\,\left(\int_{0}^{u_{t}-1}\,\rho_{\frac{1}{2}}(z)\omega(u_{t}-z )\,dz+\rho_{\frac{1}{2}}(u_{t})\right)\; \times\\[0.5cm]
&&\left(1+O\left(\frac{\log(\frac{\log t}{\log y}+1)}{\log y}+\frac{1}{\sqrt{\log y}}\right)\;\right)\;\frac{dt}{t}\cdot
\end{array}$$
A change of variable  $s=\frac{\log t}{\log y}$ gives

$$
\begin{array}{lcl}
R_{2,1}&=&-x \int_{0}^{uv}  \rho_{\frac{1}{2}}(s)\left(\int_{0}^{u-s-1}\rho_{\frac{1}{2}}(z)\omega(u-s-z )dz+\rho_{\frac{1}{2}}(u-s)\right)\; \times\\[0.5cm]
&&\left(1+O\left(\frac{\log (s+1)}{\log y}+\frac{1}{\sqrt{\log y}}\right)\right)ds\\[0.5cm]
&=&- x \,\int_{0}^{uv}  \rho_{\frac{1}{2}}(s)\,\left(\int_{0}^{u-s-1}\,\rho_{\frac{1}{2}}(z)\omega(u-s-z )\,dz\right)\,ds-x\,\,\int_{0}^{uv}  \rho_{\frac{1}{2}}(s)\rho_{\frac{1}{2}}(u-s) ds\\[0.5cm]
& +&O\left(\frac{x\,\log(u+1)}{\log y}+\frac{x}{\sqrt{\log y}}\right)\\[0.5cm]
\end{array}
$$
Now consider $R_{2,2}$. an easy calculation gives 

$$D_{2}\ll \frac{1}{t^{2} (u_{t}-1)\log y}$$

 $$
 \begin{array}{lcl}
 D_{2}&=&-\frac{1}{t^{2} (u_{t}-1)\log y}\left(\int_{0}^{u_{t}-1}\rho_{\frac{1}{2}}(z)\,\omega(u_{t}-z )\, dz\right)\;dt\\[0.5cm]
 &-&\frac{1}{t^{2}(u_{t}-1) \log y}\left(\int_{0}^{u_{t}-1}z\rho'_{\frac{1}{2}}(z)\,\omega(u_{t}-z )\, dz\right)\;dt\\[0.5cm]
 &-&\frac{1}{t^{2} (u_{t}-1)\log y}\left(\int_{0}^{u_{t}-1}(u_{t}-z)\rho_{\frac{1}{2}}(z)\,\omega'(u_{t}-z )\, dz\right)\;dt\\[0.5cm]
 &-&\frac{1}{t^{2} \log y}\,\rho'_{\frac{1}{2}}(u_{t})\\[0.5cm]
 &:=& D_{2,1}+D_{2,2}+D_{2,3}+D_{2,4}.
  \end{array}$$
We have
$$R_{2,2}=\frac{x}{\sqrt{\log y}}\,B(1) \;\int_{1}^{x^{v}}  \left(\sum_{n\leqslant t, P(n)\leqslant y}g(n)\right) \left(D_{2,1}+D_{2,2}+D_{2,3}+D_{2,4}\right)\;dt=\widetilde{D}_{1}+\widetilde{D}_{2}+\widetilde{D}_{3}+\widetilde{D}_{4}$$
where 
$$\widetilde{D}_{i}=\frac{x}{\sqrt{\log y}}\,B(1) \;\int_{1}^{x^{v}}  \left(\sum_{n\leqslant t, P(n)\leqslant y}g(n)\right) D_{2,i}$$
 for $i\in \{1,2,3,4\}$.
By using the same method as previously, we get
$$
\widetilde{D}_{1}\ll\frac{x}{\log y}\int_{0}^{uv}\frac{\rho_{\frac{1}{2}}(s)}{u-s-1}\,\left(\int_{0}^{u-s-1}\,\rho_{\frac{1}{2}}(z)\omega(u-s-z )\,dz\right)\,ds.$$
For  $s\geqslant u-2$ we have   $\omega (u-s-z)=0$ since  $0 \leqslant z \leqslant u-s-1\leqslant1$ and for  $s\leqslant u-2$, we have $u-s-1\geqslant1$. It follows  that
$$\widetilde{D}_{1}\ll \frac{x}{\log y}\int_{0}^{u-2}\rho_{\frac{1}{2}}(s)ds\,\int_{0}^{u-1}\,\rho_{\frac{1}{2}}(z)dz
\ll\frac{x}{\log y}\cdot$$
The study of  $\widetilde{D}_{2}$ , $\widetilde{D}_{3}$ and $\widetilde{D}_{4}$ is similar by using propreties of  $\rho'_{\frac{1}{2 }}$ and $\omega'$. We get the same result. We omit details.
Assembling these estimates we get our result in the considered range. 
In the complementary range  $x/d< y\leqslant x$, proof is similar, we omit it. 
 \hfill{
$\square$}
\section{Proof of  Theorem 2.2}
\subsection{Lemmas}
We need a weak form of the following lemma
\begin{Lem}
For each fixed integer $N\geqslant1$ and uniformly for $w>1$ we have
$$\int_{w}^{\infty } \rho_{\frac{1}{2}}(t) dt=2\sum_{k=1}^{N}(w+k)\rho_{\frac{1}{2}}(w+N)+ O\left(\frac{\rho_{\frac{1}{2}}(w)}{w^{N}}\right).$$
\end{Lem}
\textbf{Proof} Set
$$I(w):=\int_{w}^{\infty } \rho_{\frac{1}{2}}(t) dt.$$
By change of variable  $ z=t+1$ and by using  the differerential equation satisfied by  $\rho_{\frac{1}{2}}$ we obtain
\begin{eqnarray}
I(w)=\int_{w+1}^{\infty } \rho_{\frac{1}{2}}(z-1) dz= -2\int_{w+1}^{\infty } z\rho'_{\frac{1}{2}}(z) dz-\int_{w+1}^{\infty } \rho_{\frac{1}{2}}(z) dz  = -2 J(w)-I(w+1).
\end{eqnarray}
with  $$J(w):=\int_{w+1}^{\infty } z\rho'_{\frac{1}{2}}(z) dz.$$
Integration by parts gives
$$J(w)=-(w+1)\rho_{\frac{1}{2}}(w+1)-\int_{w+1}^{\infty } \rho_{\frac{1}{2}}(z)dz=-(w+1)\rho_{\frac{1}{2}}(w+1)-I(w+1).$$
Substituting in  (11) we obtain
\begin{eqnarray}
I(w)-I(w+1)=2(w+1)\rho_{\frac{1}{2}}(w+1).
\end{eqnarray}
And by iteration we obtain for every integer  $k\geqslant1$,
$$I(w+k-1)-I(w+k)=2(w+k)\rho_{\frac{1}{2}}(w+k),$$
Summing these inequalities, we obtain 
$$I(w)=2(w+1)\rho_{\frac{1}{2}}(w+1)+2\sum_{k=2}^{\infty}(w+k)\rho_{\frac{1}{2}}(w+k),$$

The lemma follows from this formula and the result
$$\rho_{\frac{1}{2}}(w+k)=O\left(\frac{\rho_{\frac{1}{2}}(w)}{w^{k}}\right),$$
 uniformly for $w>1$ and  $k>0$
which is a consequence of the asymptotic formula for  $\rho_{\frac{1}{2}}(w)$ (see  [8]). 
 \hfill{
$\square$}
\begin{Lem}
For $v\in[0,1]$ and as  $u\rightarrow \infty$ we uniformly have  
$$H(u,v):= \frac{1}{e^{\gamma}}  \int_{0}^{uv}\left( \int_{0}^{u-s-1}\rho_{\frac{1}{2}}(w)dw\right)\,\rho_{\frac{1}{2}}(s)ds= \frac{1}{\sqrt{e^{\gamma}}}  \int_{0}^{uv}\rho_{\frac{1}{2}}(w)dw+O\left( u\rho(u)\right),$$
where  $\gamma$ is Euler's constant.
\end{Lem}
\textbf{Proof }
 We have 
$$
\begin{array}{lcl}
H(u,v)&=& \frac{1}{e^{\gamma}}  \int_{0}^{uv}\left( \int_{0}^{+\infty} \rho_{\frac{1}{2}}(w)dw -\int_{u-s-1}^{\infty}\rho_{\frac{1}{2}}(w)dw\right)\,\rho_{\frac{1}{2}}(s)ds\\[0.5cm]
&=& \frac{1}{\sqrt{e^{\gamma}}}  \int_{0}^{uv} \rho_{\frac{1}{2}}(s)ds - \frac{1}{e^{\gamma}}  \int_{0}^{uv}\left(\int_{u-s-1}^{\infty}\rho_{\frac{1}{2}}(w)dw\right)\,\rho_{\frac{1}{2}}(s)ds,$$
\end{array}$$
since by using Laplace transform, we have  $$\int_{0}^{+\infty} \rho_{\frac{1}{2}}(w)dw=\widehat{\rho_{\frac{1}{2}}}(0)=(\widehat{\rho}(0))^{1/2}=\sqrt{e^{\gamma}}.$$
We apply Lemma 4.1 in a weak version we obtain
$$
\begin{array}{lcl}
H(u,v)&=& \frac{1}{\sqrt{e^{\gamma}}}  \int_{0}^{uv} \rho(s)ds+O\left( u \int_{0}^{uv} \rho_{\frac{1}{2}}(s)\,\rho_{\frac{1}{2}}(u-s)\,ds \right)\\
&=&\frac{1}{\sqrt{e^{\gamma}}}  \int_{0}^{uv} \rho(s)ds+O(u\rho(u)),
\end{array}
$$
since
$$ \int_{0}^{uv} \rho_{\frac{1}{2}}(s)\,\rho_{\frac{1}{2}}(u-s)\,ds\leqslant  \int_{0}^{u} \rho_{\frac{1}{2}}(s)\,\rho_{\frac{1}{2}}(u-s)\,ds =(\rho_{\frac{1}{2}}\ast \rho_{\frac{1}{2}})(u)\ll\rho(u).$$
\hfill{
$\square$}\\[0.5cm]
\subsection{Proof of Theorem 2.2}
We have 
$$F(u,v)=\int_{0}^{uv} \left(\int_{0}^{u-s-1}\rho_{\frac{1}{2}}(z)\omega(u-s-z )\,dz\right)\rho_{\frac{1}{2}}(s)\,ds+O(\rho(u)).$$
 And
 $$
 \begin{array}{lcl}
\int_{0}^{uv} \left(\int_{0}^{u-s-1}\rho_{\frac{1}{2}}(z)\omega(u-s-z )\,dz\right)\rho_{\frac{1}{2}}(s)\,ds= \frac{1}{\sqrt{e^{\gamma}}}  \int_{0}^{uv}\rho_{\frac{1}{2}}(s)ds-\frac{1}{\sqrt{e^{\gamma}}}  \int_{0}^{uv}\rho_{\frac{1}{2}}(s)\left(\int_{u-s-1}^{\infty}\rho_{\frac{1}{2}}(z)dz\right) ds+
  O\left(  \int_{0}^{uv} \rho_{\frac{1}{2}}(s) \left(\int_{0}^{u-s-1}\rho_{\frac{1}{2}}(z)\rho(u-s-z)dz\right)\,ds\   \right)
   \end{array}
   $$
 By using  formula  (3) and Lemma  4.2, we get
 $$
 \begin{array}{lcl}
&& \int_{0}^{uv} \left(\int_{0}^{u-s-1}\rho_{\frac{1}{2}}(z)\omega(u-s-z )\,dz\right)\rho_{\frac{1}{2}}(s)\,ds=
\frac{1}{e^{\gamma}}\int_{0}^{uv} \rho_{\frac{1}{2}}(s) \left(\int_{0}^{u-s-1}\rho_{\frac{1}{2}}(z)dz\right)\,ds\, \\[0.5cm]
 &+& O\left(  \int_{0}^{uv} \rho_{\frac{1}{2}}(s) \left(\int_{0}^{u-s-1}\rho_{\frac{1}{2}}(z)\rho(u-s-z)dz\right)\,ds\   \right)\\[0.5cm]
&=&\frac{1}{\sqrt{e^{\gamma}}}  \int_{0}^{uv}\rho_{\frac{1}{2}}(s)ds+O(u\rho(u))+O\left( \rho_{2}(u)\right)\\[0.5cm]
&=&\frac{1}{\sqrt{e^{\gamma}}}  \int_{0}^{uv}\rho_{\frac{1}{2}}(s)ds+O\left( \rho_{2}(u)\right)
\end{array}
$$
since
$$
\begin{array}{lcl}
\int_{0}^{uv} \rho_{\frac{1}{2}}(s) \left(\int_{0}^{u-s-1}\rho_{\frac{1}{2}}(z)\rho(u-s-z)dz\right)\,ds& \leqslant &\int_{0}^{uv} \rho_{\frac{1}{2}}(s)(\rho_{\frac{1}{2}}\ast \rho)(u-s)\,ds\\
&\leqslant&(\rho_{\frac{1}{2}}\ast \rho_{\frac{1}{2}}\ast \rho)(u)=(\rho\ast \rho)(u)=\rho_{2}(u),
\end{array}
$$
and   $u\rho(u)\ll \rho_{2}(u).$
\hfill{
$\square$}
\section{Proof of Theorem 2.3}
\subsection{Lemmas}
\begin{Lem}  For $0\leqslant \xi <1$ we have
$$R(\xi):=\frac{1}{\pi} \int_{0}^{\xi} \int_{0}^{-s+\xi} \frac{ds'\,dz'} {\sqrt{s'}\sqrt{z'}(1-s'-z')}=
-\log (1-\xi)\cdot$$
\end{Lem}
\textbf{Proof}
We use change of variables $(s',z')\mapsto (w, r)=(\frac{s'}{s'+z'},s'+z')$. Then we have  a $dsdz=rdwdr$.
We obtain
 $$
R(\xi) =\left(\frac{1}{\pi}\int_{0}^{1}w ^{-\frac{1}{2}} (1-w)^{-\frac{1}{2}}\;dw\right)\left(  \int_{0}^{\xi}(1-r)^{-1} \;dr \right)=-\log (1-\xi),
$$
since
$$\frac{1}{\pi}\int_{0}^{1}w ^{-\frac{1}{2}} (1-w)^{-\frac{1}{2}}\;dw=\frac{B(\frac{1}{2},\frac{1}{2})}{\pi}=\frac{\Gamma^{2}(\frac{1}{2})}{\pi}=1.$$
\hfill{
$\square$}
\begin{Lem} For $\beta \leqslant 1-\frac{1}{u}$ we have
$$
\begin{array}{lcl}
S(\beta)&:=&\frac{1}{\pi}\int_{0}^{\beta} \left(\int_{0}^{-s+1-\frac{1}{u}}\frac{ds \,dz}{\sqrt{s}\,\sqrt{z}\,(1-s-z)}\right)\\[0.5cm]
&=&\frac{2}{\pi}\left(\;\log(u)+\log(1-\beta)\;\right)\arcsin\left(  \sqrt{\frac{u\beta}{u-1}} \,  \right)-\log(1-\beta)\cdot
\end{array}
$$
\end{Lem}
\textbf{Proof} We write $S(\beta)=I_{1}-I_{2}$ with
$$I_{1}=\frac{1}{\pi}\int_{0}^{1-\frac{1}{u}} \left(\int_{0}^{-s+1-\frac{1}{u}}\frac{ds \,dz}{\sqrt{s}\,\sqrt{z}\,(1-s-z)}\right),\,I_{2}=\frac{1}{\pi}\int_{\beta}^{1-\frac{1}{u}} \left(\int_{0}^{-s+1-\frac{1}{u}}\frac{ds \,dz}{\sqrt{s}\,\sqrt{z}\,(1-s-z)}\right)\cdot$$
From Lemma 5.1, $I_{1}=R(1-\frac{1}{u})=\log(u)$. Let us study  $I_{2}$. By change of  variables $(s',z')\mapsto(s,z)=(s'+\beta,z')$ we obtain
$$I_{2}=\frac{1}{\pi}\int_{0}^{1-\frac{1}{u} -\beta} \left(\int_{0}^{-s'+1-\frac{1}{u} -\beta}\frac{ds' \,dz'}{\sqrt{s'+\beta}\,\sqrt{z'}\,(1-s'-\beta-z')}\right)\cdot$$
Now we put the change of variables
$(s'+\beta)=rw$ et $s'+\beta+z'=r$ then we have     $ds'dz' =rdwdr$ and 
$$I_{2}=\frac{1}{\pi}\int_{\frac{u\beta}{u-1}}^{1}  \frac{dw}{\sqrt{w} \sqrt{1-w} }\int_{\beta}^{1-\frac{1}{u}}\frac{dr}{1-r}\cdot$$
Since
$$\int_{\beta}^{1-\frac{1}{u}}\frac{dr}{1-r}=\log(u)+\log(1-\beta),$$
then by change of variable $t=\sqrt{w}$, we get
$$\frac{1}{\pi}\int_{\frac{u\beta}{u-1}}^{1}  \frac{dw}{\sqrt{w} \sqrt{1-w}} =\frac{2}{\pi}\int_{\sqrt{\frac{u\beta}{u-1}}}^{1} \frac{dt}{ \sqrt{1-t^{2}} }=1-\frac{2}{\pi}\arcsin\left(  \sqrt{\frac{u\beta}{u-1}} \,  \right)\cdot$$
Finally, 
$$I_{2}=\left(\;\log(u) +\log(1-\beta)\;\right)\left(1-\frac{2}{\pi}\arcsin \left(  \sqrt{\frac{u\beta}{u-1}} \,\right) \right)\cdot$$
It follows
$$S(\beta)=\frac{2}{\pi}( \log(u)+\log(1-\beta))\arcsin \left(  \sqrt{\frac{u\beta}{u-1}}\right) -\log(1-\beta).$$
\hfill{
$\square$}
\begin{Lem} 
For $0<w\leqslant1$ we have  $$\rho_{\frac{1}{2}}(w)=\frac{1}{\sqrt{\pi} \sqrt{w}},$$
and for  $1\leqslant w\leqslant2$, we have 
$$\rho_{\frac{1}{2}}(w)=\frac{1}{\sqrt{\pi} \sqrt{w}} 
- \frac{\log(\sqrt{w}+\sqrt{w-1})}{\sqrt{\pi} \sqrt{w}}.$$
\end{Lem}
\textbf{Proof} The first formula is the definition of $\rho_{\frac{1}{2}}$  for $0< w \leqslant1$ and the second one follows  from differential equation
satisfied by  $\rho_{\frac{1}{2}}$ pour  $1\leqslant w \leqslant2 $. \hfill{
$\square$}\\[0.2cm]
Let us consider the integral
$$I:=\int_{0}^{uv} \left(\int_{0}^{u-s-1}\rho_{\frac{1}{2}}(z)\omega(u-s-z )\,dz\right)\rho_{\frac{1}{2}}(s)\,ds,$$
 for  $1\leqslant u\leqslant 2$. We notice  that on one hand  
 for  $s>u-1$ 
$$\int_{0}^{u-s-1}\rho_{\frac{1}{2}}(z)\omega(u-s-z )\,dz=0,$$
since $\rho_{\frac{1}{2}}(z)=0$ pour $z\leqslant0$ and then we can restrict the study of the l'integral on  $s$ at the interval $[0,M]$ where $M=\min\{u-1,uv\}\leqslant 1$.
And on the other hand, $0\leqslant z\leqslant u-s-1$ then $1\leqslant u-s-z\leqslant u-s\leqslant u\leqslant 2$  and   $\omega (u-s-z)=1/(u-s-z)$.  It follows
 \begin{eqnarray}
I=\frac{1}{\pi}\int_{0}^{M} \left(\int_{0}^{-s+u-1}\frac{dz}{\sqrt{z}\,(u-s-z)}\right)\frac{ds}{\sqrt{s}}\cdot
\end{eqnarray}
 We will give two expessions of  $I$.
By the change of variable $t=\sqrt{z}$ in the inner integral $(13)$ we get
\begin{eqnarray}
I=\frac{2}{\pi}\,\int_{0}^{M}\frac{\log(\sqrt{u-s}+\sqrt{u-s-1})}{\sqrt{s} \sqrt{u-s}}\,ds,
\end{eqnarray}
and by the change of variables $(s',z')\mapsto (s,z)=(us',uz')$ and by puting  $M'=\frac{M}{u}=\min \{ 1-\frac{1}{u},v \}$
in  (13), we obtain
\begin{eqnarray}
 I=\frac{1}{\pi}\int_{0}^{M'}\int_{0}^{-s'+1-\frac{1}{u}}\frac{ds'\,dz'}{\sqrt{s'}  \sqrt{z'}(1-s'-z')}\cdot
\end{eqnarray}
By using notations of  lemma 5.2 and the expressions  (13) et (14) of $I$, we obtain
\begin{eqnarray}
I=
\left \lbrace
\begin{array}{lcl}
S(1-\frac{1}{u})\quad \text{si}\quad M=u-1,\\[0.2cm]
S(v)\qquad \quad \, \text{si}\quad M=uv.
\end{array}
\right.
\end{eqnarray}
Now we will study the integral
$$\int_{0}^{uv} \rho_{\frac{1}{2}}(s)\,\rho_{\frac{1}{2}}(u-s)\,ds,\quad (1\leqslant u \leqslant 2).$$

\begin{Lem}
For $v\in[0,\frac{u-1}{u}]$ with  $1\leqslant u \leqslant2$ we have
$$\int_{0}^{uv} \rho_{\frac{1}{2}}(s)\,\rho_{\frac{1}{2}}(u-s)\,ds=\frac{2}{\pi}\arcsin(\sqrt{v})-\frac{1}{2} S(v).$$
\end{Lem}
\textbf{Proof}
For $v\in[0,\frac{u-1}{u}]$ we have  $uv\leqslant u-1$ then $0\leqslant s \leqslant uv\leqslant u-1\leqslant1$ and
$1\leqslant u-s \leqslant 2$. By using Lemma 5.3 we get
$$\int_{0}^{uv} \rho_{\frac{1}{2}}(s)\,\rho_{\frac{1}{2}}(u-s)\,ds=\frac{1}{\pi}\int_{0}^{uv} \frac{ds}{\sqrt{s}\sqrt{u-s}}-\frac{1}{\pi}\int_{0}^{uv} \frac{\log(\sqrt{u-s}+\sqrt{u-s-1})}{\sqrt{s}\sqrt{u-s}}\,ds:=J_{1}-J_{2}.$$
By the change of variable $t=\sqrt{s}$ we have
$$J_{1}=\frac{2}{\pi}\int_{0}^{\sqrt{v}} \frac{dt}{\sqrt{1-t^{2}}}=\frac{2}{\pi}\arcsin(\sqrt{v}).$$
From (16) we have $J_{2}=\frac{1}{2}S(v)$. This completes the proof of Lemma 5.4.
\hfill{
$\square$}
\begin{Lem}
For $v\in[\frac{u-1}{u},\frac{1}{u}]$ with $1\leqslant u\leqslant 2$ we have
$$\int_{0}^{uv} \rho_{\frac{1}{2}}(s)\,\rho_{\frac{1}{2}}(u-s)\,ds=\frac{2}{\pi}\arcsin(\sqrt{v})-\frac{1}{2}S(1-\frac{1}{u}).$$
\end{Lem}
\textbf{Proof}  For  $v\in[\frac{u-1}{u},\frac{1}{u}]$we have  $u-1 \leqslant uv \leqslant1$. We write
$$\int_{0}^{uv} \rho_{\frac{1}{2}}(s)\,\rho_{\frac{1}{2}}(u-s)\,ds=\int_{0}^{u-1} \rho_{\frac{1}{2}}(s)\,\rho_{\frac{1}{2}}(u-s)\,ds+\int_{u-1}^{uv} \rho_{\frac{1}{2}}(s)\,\rho_{\frac{1}{2}}(u-s)\,ds=J_{1}+J_{2}\cdot$$
Consider $J_{1}$. As in previous Lemma we have
$$
\begin{array}{lcl}
J_{1}&=&\frac{1}{\pi}\int_{0}^{u-1} \frac{ds}{\sqrt{s}\sqrt{u-s}}-\frac{1}{\pi}\int_{0}^{u-1} \frac{\log(\sqrt{u-s}+\sqrt{u-s-1})}{\sqrt{s}\sqrt{u-s}}\,ds\\[0.5cm]
&=& \frac{2}{\pi} \arcsin(\sqrt{\frac{u-1}{u}})-\frac{1}{2}S(1-\frac{1}{u}).
\end{array}
$$
Now consider  $J_{2}$. We have  $u-1\leqslant s\leqslant uv\leqslant1$ et $u(1-v)\leqslant u-s\leqslant 1$. It follows
$$J_{2}=\frac{1}{\pi}\int_{u-1}^{uv} \frac{ds}{\sqrt{s}\sqrt{u-s}}=\frac{2}{\pi}\arcsin(\sqrt{v})-\frac{2}{\pi}\arcsin(\sqrt{\frac{u-1}{u}}).$$
This completes the proof of Lemma 5.5.
\hfill{
$\square$}
\begin{Lem}
For $v\in[\frac{1}{u}, 1]$ with $1\leqslant u\leqslant 2$ we have 
$$\int_{0}^{uv} \rho_{\frac{1}{2}}(s)\,\rho_{\frac{1}{2}}(u-s)\,ds=-\frac{3}{2}S(1-\frac{1}{u})+\frac{2}{\pi}\arcsin(\sqrt{v})+S(1-v).$$
\end{Lem}
 \textbf{Proof}
As $u\leqslant 2$ we have $\frac{u-1}{u}\leqslant \frac{1}{u}\leqslant v\leqslant1$ and then  $u-1\leqslant uv\leqslant u$. We write
$$ \int_{0}^{uv} \rho_{\frac{1}{2}}(s)\,\rho_{\frac{1}{2}}(u-s)\,ds=\int_{0}^{u-1} \rho_{\frac{1}{2}}(s)\,\rho_{\frac{1}{2}}(u-s)\,ds+\int_{u-1}^{uv} \rho_{\frac{1}{2}}(s)\,\rho_{\frac{1}{2}}(u-s)\,ds=J_{1}+J_{2}\cdot$$
$J_{1}$ has been studied in lemma 5.5. It remains to calculate $J_{2}$. As  $\frac{1}{u} \leqslant v$ then $uv\geqslant1$. We write
$$J_{2}=\int_{u-1}^{1} \rho_{\frac{1}{2}}(s)\,\rho_{\frac{1}{2}}(u-s)\,ds+\int_{1}^{uv} \rho_{\frac{1}{2}}(s)\,\rho_{\frac{1}{2}}(u-s)\,ds := J_{2,1}+J_{2,2}.$$
We have
$$J_{2,1}=\frac{ 1}{\pi}\int_{u-1}^{1}\frac{ds}{\sqrt{s}\sqrt{u-s}}=\frac{2}{\pi}\arcsin(\frac{1}{\sqrt{u}})-\frac{2}{\pi}\arcsin(\sqrt{\frac{u-1}{u}}).$$
 We write
$$J_{2,2}=\frac{1}{\pi}\int_{1}^{uv}\frac{ds}{\sqrt{s}\sqrt{u-1}}-\frac{1}{\pi}\int_{1}^{uv}\frac{\log(\sqrt{s}+\sqrt{s-1})}{\sqrt{s}\sqrt{u-s}}\,ds:=\widehat{J_{1}}-\widehat{J_{2}}.$$
We have
$$\widehat{J_{1}}=\frac{2}{\pi}\arcsin(\sqrt{v})-\frac{2}{\pi}\arcsin(\frac{1}{\sqrt{u}}).$$
Now it remains to study $\widehat{J_{2}}$. We put the change of variable  $s'=u-s$. We get
$$\widehat{J_{2}}=\frac{1}{\pi}\int_{u(1-v)}^{u-1}\frac{\log(\sqrt{u-s'}+\sqrt{u-s'-1})}{\sqrt{s'}\sqrt{u-s'}}\,ds'$$
As $v\geqslant \frac{1}{u}$ we have $u(1-v)\leqslant u-1\leqslant 1$. Therefore using notations in (14) (15) and (16) and Lemmas 5.1 et 5.2 we obtain
$$\widehat{J_{2}}=S(1-\frac{1}{u})-S(1-v).$$
We complete the proof by grouping  different  above estimates.
\hfill{
$\square$}\\[0.5cm]
\subsection{Proof of Theorem 2.3}
Theorem 2.3 follows  from  (15) and different Lemmas  of section 5.   \hfill{
$\square$}\\[0.5cm]
\selectlanguage{english}

\flushleft{Université du Burundi \;\hfill Université de Limoges\;\quad \quad \qquad  \,\, \\
\;Faculté des sciences \hfill UMR-CNRS 7252 \qquad \qquad \qquad  \qquad \qquad\\
Avenue de l'UNESCO n$^{\mbox{\tiny 0 }}1$\quad\hfill \;123 avenue Albert Thomas \qquad  \qquad \qquad\\
B.P 2700 Bujumbura, Burundi \; \hfill 87060 Limoges Cedex, France \,\; }
\end{document}